\documentclass[12pt]{article}
\usepackage{pifont}

\usepackage{amsfonts}
\usepackage{latexsym}
\usepackage{amsmath}
\usepackage{amssymb}
\usepackage{color}

\def\o{\Omega}
\def\z{\zeta}
\def\a{\alpha}
\def\ep{\varepsilon}
\def\p{\varphi}
\def\d{\delta}

\def\g{\gamma}

\title{ Viscosity solutions to second order parabolic PDEs on
Riemannian manifolds }

\date{June 7, 2010}

\author{ Xuehong Zhu}

\begin{document}

\maketitle

\begin{abstract}
In this work we consider viscosity solutions to second order
  parabolic PDEs $u_{t}+F(t,x,u,du,d^{2}u)=0$ defined on compact Riemannian manifolds with boundary
  conditions. We prove comparison, uniqueness and existence results
  for the solutions. Under the
  assumption that the manifold $M$ has nonnegative sectional
  curvature, we get the finest results. If one additionally requires
  $F$ to depend on $d^{2}u$ in a uniformly continuous manner, the
  assumptions on curvature can be thrown away.\\
\par  $\textit{Keywords: Second order parabolic PDEs;\\ Riemannian manifold;
Viscosity solution.}$
\end{abstract}

%% \linenumbers

%% main text

\section{Introduction}\label{sec:intro}

\qquad Since the theory of viscosity solutions was introduced by M.
G. Crandall and P. L. Lions in the 1980's, it has been found that
it's a very natural concept for a generalization of classical
solutions. This theory has been applied widely and was enriched and
expanded by many mathematicians. We would refer the reader to [1]
and the references therein.

It should be natural to generalize the theory to problems on
Riemannian manifolds since many functions arise from geometrical
problems. But little is known about this field. Until recently, D.
Azagra, J. Ferrera and B. Sanz[2] gave a work about Dirichlet
problem on a complete Riemannian manifold with some restrictions on
curvature. Almost at the same time, a few results about parabolic
PDEs on Riemannian manifolds without boundary conditions were given
in the appendix of [3].

We consider the following Cauchy-Dirichlet problem of the form:

$$\left\{
\begin{array}{ll}
(E)&   u_{t}+F(t,x,u,du,d^{2}u)=0\mbox{ \ in \ }(0,T)\times \o\\
(BC) & u(t,x)=h(t,x), (t,x)\in [0,T)\times \partial \o, \\
(IC) & u(0,x)=\psi(x),  x\in \bar{\o},
\end{array}
\right.
 \eqno{(1.1)}
$$
where $u$ is a function of $(t,x)$: $[0,T)\times M\rightarrow R$ and
$M$ is finite-dimensional complete Riemannian manifold. $du, d^{2}u$
mean $d_{x}u(t,x)$ and $d^{2}_{x}u(t,x)$. $\o$ in $M$ is open and
bounded. $T>0$, $h\in C([0,T)\times \bar{\o})$ and $\psi\in
C(\bar{\o})$ are given.

Since there are not many intrinsic differences between elliptic and
parabolic PDEs on manifolds, through an important theorem which is a
Riemannian version of a result in Euclidean space in [1], we can get
several comparison and uniqueness results which are parabolic
versions of that in [2].

When we are faced with the existence of viscosity solutions,
Perron's method is a good choice. Detailed research has been done in
[1] about Dirichlet problem in Euclidean space. We will generalize
two lemmas in [1] to parabolic PDEs on  Riemannian manifolds and
then get our existence result.

The next section is about the properties of second order parabolic
viscosity subdifferentials on Riemannian manifolds and an important
theorem for the proofs of comparison results; Section 3 are devoted
to comparison results of our PDEs (1.1); Finally we prove the
existence result of (1.1).

\section{Second order parabolic viscosity
subdifferentials on Riemmanian manifolds}\label{sec:intro}

{\bf Definition 2.1.} Let $M$ be a finite-dimensional Riemannian
manifold, and $f : (0,T)\times M \rightarrow R$ a lower
semicontinuous function. We define the second order parabolic subjet
 of $f$ at a point $(t,x)\in (0,T)\times M$ by
$$
\begin{array}{ll} & {\cal{P}}^{2,-}f(t,x)\\
=&\{(\p_{t}(t,x),d\p(t,x),d^{2}\p(t,x)) : \p\in C^{1,2}([0,T]\times
M),\\
&\mbox{ \ \ \ \ \ \ \ \ \ \ \ \ \ \ \ \ \ \ \ \ \ \ \ \ \ \ \ \ \ \
 \ \ \ \ \ \ } f-\p \mbox{ \ attains a local minimum 0 at \ } (t,x)\}.
\end{array}
$$
If $(p,\z,A) \in {\cal{P}}^{2,-}f(t,x)$, we will say that at the
point $(t,x)$, $p$ is a first order subdifferential of $f$ w.r.t.
$t$, $\z$ is a first order subdifferential of $f$ w.r.t. $x$ and $A$
is a second order subdifferential of $f$ w.r.t. $x$.

Similarly, for an upper semicontinuous function $f : (0,T)\times M
\rightarrow R$, we define the second order parabolic superjet
 of $f$ at a point $(t,x)$ by
$$
\begin{array}{ll}
&{\cal{P}}^{2,+}f(t,x)\\
=&\{(\p_{t}(t,x),d\p(t,x),d^{2}\p(t,x)) : \p\in C^{1,2}([0,T]\times
M),\\
&\mbox{ \ \ \ \ \ \ \ \ \ \ \ \ \ \ \ \ \ \ \ \ \ \ \ \ \ \ \ \ \ \
 \ \ \ \ \ \ } f-\p \mbox{ \ attains a local maximum 0 at \ } (t,x)\}.
\end{array}$$
 Observe that ${\cal{P}}^{2,-}f(t,x)$ and $
{\cal{P}}^{2,+}f(t,x)$ are subsets of $R\times TM_{x}^{*}\times
{\cal{L}}_{s}^{2}(TM_{x})$, where $TM_{x}^{*}$ stands for the
cotangent space of $M$ at a point $x$, $TM_{x}$ stands for the
tangent space at $x$ and ${\cal{L}}_{s}^{2}(TM_{x})$ denotes the
symmetric bilinear forms on $TM_{x}$. It is also clear that
${\cal{P}}^{2,-}f(t,x)=-{\cal{P}}^{2,+}(-f)(t,x)$.

In the sequel $M$ will always denotes an $n$-dimensional Riemannian
manifold. We need the following several results for subjets which
also hold, with obvious modification, for superjets.

{\bf Proposition 2.2.} Let $f : (0,T)\times M\rightarrow R$ be a
lower semicontinuous function. Let $p\in R, \z\in TM_{x}^{*},A\in
{\cal{L}}_{s}^{2}(TM_{x}), (t,x)\in (0,T)\times M$. The following
statements are equivalent:
$$ \begin{array}{l}
\mbox{(i)} (p,\z,A)\in {\cal{P}}^{2,-}f(t,x).\\
\mbox{(ii)} f(s,\exp_{x}(v))\geq
f(t,x)+p(s-t)+\langle\z,v\rangle_{x}+\frac{1}{2}\langle
Av,v\rangle_{x}+o(|s-t|+\|v\|^{2}), \\
\mbox{ \ \ \ \ \ \ \ \ \ \ \ \ \ \ \ \ \ \ \ \ \ \ \ \ \ \ \ \ \ \
 \ \ \ \ \ \ }\mbox{ \ \ \ \ \ \ \ \ \ \ \ \ \ \ \ \ \ \ \ \ \ \ \ \ \ \ \ \ \ \
 \ }
\mbox{as \ } |s-t|\rightarrow 0, \|v\|\rightarrow 0.
\end{array} $$

{\bf Remark 2.3.} This result is mainly due to [2] with a slight
modification: $f$ is also a function of time $t$. Note that if
$\p(t,x): (0,T)\times M\rightarrow R$ is differentialble w.r.t $t$
in $(0,T)$. Set $\psi(t,v)=\p(t,\exp_{x}(v))$, where $v$ is valued
on a neighborhood of $0_{x}$ in $TM_{x}$.  Then
$\p_{t}(s,y)=\psi_{t}(s,\exp_{x}^{-1}(y))$. So obviously we can use
almost the same method as Proposition 2.2 in [2] to get the proof of
our proposition. The following Corollary 2.4 and Proposition 2.6 are
also analogues of that in [2]. So we omit the proofs.

{\bf Proposition 2.4.} Let $f : (0,T)\times M\rightarrow R$ be a
lower semicontinuous function, and consider $p\in R, \z\in
TM_{x}^{*},A\in {\cal{L}}_{s}^{2}(TM_{x}), (t,x)\in (0,T)\times M$.
Set $\bar{f}(t,v)=f(t,\exp_{x}(v))$. Then
$$
(p,\z,A)\in {\cal{P}}^{2,-}f(t,x)\Leftrightarrow (p,\z,A) \in
{\cal{P}}^{2,-}\bar{f}(t,0_{x}).
$$

{\bf Definition 2.5. }Let $f : (0,T)\times M\rightarrow R$ be a
lower semicontinuous function, and $(t,x)\in (0,T)\times M$. We
define
$$
\begin{array}{ll}
\bar{{\cal{P}}}^{2,-}f(t,x)=& \{(p,\z,A)\in R\times TM_{x}^{*}\times
{\cal{L}}_{s}^{2}(TM_{x}) :\\
& \mbox{ \ }\exists (t_{n},x_{n},p_{n},\z_{n},A_{n})\in (0,T)\times
M\times
R\times TM_{x_{n}}^{*}\times {\cal{L}}_{s}^{2}(TM_{x_{n}}), \\
& \mbox{ \ }s.t. (p_{n},\z_{n},A_{n})\in
{\cal{P}}^{2,-}f(t_{n},x_{n}),\\
&\mbox{ \
}(t_{n},x_{n},f(t_{n},x_{n}),p_{n},\z_{n},A_{n})\rightarrow
(t,x,f(t,x),p,\z,A)\}.
\end{array}
$$
and for an upper semicontinuous function $f(t,x)$ defined on
$(0,T)\times M$ we define $\bar{{\cal{P}}}^{2,+}f(t,x)$ in the
obvious way.

{\bf Proposition 2.6. }Let $f : (0,T)\times M\rightarrow R$ be a
lower semicontinuous function, and consider $p\in R, \z\in
TM_{x}^{*},A\in {\cal{L}}_{s}^{2}(TM_{x}), (t,x)\in (0,T)\times M$.
Set $\bar{f}(t,v)=f(t,\exp_{x}(v))$. Then
$$
(p,\z,A)\in \bar{{\cal{P}}}^{2,-}f(t,x)\Leftrightarrow (p,\z,A) \in
\bar{{\cal{P}}}^{2,-}\bar{f}(t,0_{x}).
$$

The following result is the Riemannian version of Theorem 8.3 in [1]
and, as in that paper, will be the key to the proofs of comparison
and uniqueness results for viscosity solution of second order
parabolic PDEs on Riemannian manifolds. For the proof, see [3,
Theorem 3.8].

{\bf Theorem 2.7. }Let $M_{1},M_{2},...M_{k}$ be Riemannian
manifolds, and $\o_{i}\subset M_{i}$ open subsets. Define
$\o=\o_{1}\times...\times\o_{k}\subset M_{1}\times...\times
M_{k}=M$. Let $u_{i}$ be upper semicontinuous functions on
$(0,T)\times \o_{i},i=1,2,...k$; let $\p$ be in $C^{1,2}((0,T)\times
\o)$ and set
$$
w(t,x)=u_{1}(t,x_{1})+\cdots+u_{k}(t,x_{k})-\p(t,x),
$$
for $t\in (0,T), x=(x_{1},...,x_{k})\in \o$. Assume that
$(\hat{t},\hat{x})=(\hat{t},\hat{x}_{1},...,\hat{x}_{k})\in
(0,T)\times \o$ s.t.
$$
w(t,x_{1},...,x_{k})\leq
w(\hat{t},\hat{x}_{1},...,\hat{x}_{k}),\mbox{for \ } t\in (0,T),
x_{i}\in \o_{i}.
$$
Assume, moreover, that there is an $\d>0$ s.t. for every $N>0$ there
is a $C$ s.t. for $i=1,...,k$
$$
\begin{array}{l}
p_{i}\leq C \mbox{ \ whenever \ } (p_{i},\z_{i},A_{i})\in
{\cal{P}}^{2,+}u_{i}(t,x_{i}),\\
d(x_{i},\hat{x}_{i})+|t-\hat{t}|\leq \d \mbox{ \ and \ }
|u_{i}(t,x_{i})|+\|\z_{i}\|+\|A_{i}\| \leq N,
\end{array}
\eqno{(2.1)}
$$
where $d(\cdot , \cdot)$ denotes the Riemannian distance in $M$.

Then for each $\ep>0$, there are $B_{i}\in
{\cal{L}}^{2}_{s}((TM_{i})_{\hat{x}_{i}})$ such that
$$
\left\{
\begin{array}{l}
(p_{i},d_{x_{i}}\p(\hat{t},\hat{x}_{1},...,\hat{x}_{k}),B_{i})\in
\bar{{\cal{P}}}^{2,+}u_{i}(\hat{t},\hat{x}_{i}) \mbox{ \ for \ }
i=1,...,k,\\
-(\frac{1}{\ep}+\|A\|)I\leq \left(\begin{array}{ccc}
B_{1}& \cdots & 0\\
\vdots &\ddots & \vdots \\
0 & \cdots & B_{k}
\end{array}
\right)
\leq A+\ep A^{2},\\
p_{1}+\cdots+p_{k}= \p_{t}(\hat{t},\hat{x}_{1},...,\hat{x}_{k}),
\end{array}
\right.
$$
where $A=d^{2}\p(\hat{t},\hat{x})\in
{\cal{L}}_{s}^{2}(TM_{\hat{x}}).$

Now we extend the notion of viscosity solution to a parabolic
equation on a Riemannian manifold. In the sequel we will denote
$$
\chi:=\{(t,x,r,\z,A): t\in [0,T], x\in M, r\in R, \z\in TM^{*}_{x},
A\in {\cal{L}}^{2}_{s}(TM_{x})\}.
$$

{\bf Definition 2.8 (Viscosity solution). }Let $M$ be a Riemannian
manifold, and $F \in C(\chi, R)$. We say that an upper(lower)
semicontinuous function $u : [0,T)\times \bar{\o}\rightarrow R$ is a
viscosity subsolution(supersolution) of (1.1) on $[0,T)\times
\bar{\o}$, if $u(t,x)\leq(\geq) h(t,x)$ for $0\leq t<T$ and $x\in
\partial \o$ and $u(0,x)\leq(\geq) \psi(x)$ for $x\in \bar{\o}$.
Moreover, for all $(t,x)\in (0,T)\times \o$ and $(p,\z,A)\in
{\cal{P}}^{2,+}u(t,x)({\cal{P}}^{2,-}u(t,x))$,
$p+F(t,x,u(t,x),\z,A)\leq(\geq) 0$.

If $u$ is both a viscosity subsolution and a vicosity supersolution
of (1.1), we say that $u$ is a viscosity solution of (1.1).

{\bf Definition 2.9 (Degenerate ellipticity). }Let $L_{xy}$ denote
the parallel transport along the unique minimizing geodesic
connecting $x$ to $y$ (assuming $x$ is close enough to $y$ so that
this makes sense). This mapping is an isometry from $TM_{x}$ onto
$TM_{y}$ (with inverse $L_{yx}$), and it induces an isometry (which
we will still denote by $L_{xy}$), $TM^{*}_{x}\ni \z\rightarrow
L_{xy}\z\in TM^{*}_{y}$, defined by
$$
\langle L_{xy}\z,v\rangle_{y}:=\langle \z,L_{yx}v \rangle_{x}.
$$
Similarly, $L_{xy}$ induces an isometry ${\cal{L}}^{2}(TM_{x})\ni
A\rightarrow L_{xy}(A)\in {\cal{L}}^{2}(TM_{y})$ defined by
$$
\langle L_{xy}(A)v,v\rangle_{y}:=\langle A(L_{yx}v),L_{yx}v
\rangle_{x},
$$
where ${\cal{L}}^{2}(TM_{x})$ denotes the space of bilinear forms on
$TM_{x}$.

We will say that a function $F\in C(\chi,R)$ is degenerate elliptic
provided that
$$
P\leq Q \Rightarrow F(x,r,\z,Q) \leq F(x,r,\z,P),
$$
for all $x\in M, r\in R, \z\in TM_{x}^{*}, P\in
{\cal{L}}^{2}_{s}(TM_{x}), Q\in {\cal{L}}^{2}_{s}(TM_{x})$.

{\bf Definition 2.10 (Properness). }We will say that a function $F :
M\times R\times TM^{*}\times T_{2,s}(M)\rightarrow R$,
$(x,r,\z,A)\rightarrow F(x,r,\z,A)$, is proper provided

(i)$F$ is degenerate elliptic, and

(ii)$F$ is nondecreasing in the variable $r$.

\vskip 5mm

\section{ Comparison results}\label{sec:intro}
\qquad In this section and throughout the rest of the paper we will
often abbreviate saying that $u$ is an upper semicontinuous function
on a set $[0,T)\times \bar{\o}$ by writing $u\in USC([0,T)\times
\bar{\o})$. Similarly, $LSC([0,T)\times \bar{\o})$ will stand for
the set of lower semicontinuous functions on $[0,T)\times \bar{\o}$.

We denote by $i_{M}(x)$ injectivity radius of $M$ at $x$, that is
the supremum of the radius $r$ of all balls $B(0_{x},r)$ in $TM_{x}$
for which $\exp_{x}$ is a diffeomorphism from $B(0_{x},r)$ onto
$B(x,r)$. Similarly, $i(M)$ will denote the global injectivity
radius of $M$, that is $i(M)=\inf\{i_{M}(x) : x\in M\}$.

{\bf Theorem 3.1. }Let $\o$ be a bounded open subset of a complete
finite-dimensional Riemannian manifold $M$, and for each fixed $t\in
(0,T)$, $F \in C(\chi, R)$ be continuous, proper  and satisfy: there
exists a function $\omega : [0,+\infty]\rightarrow [0,+\infty]$ with
$\omega(0+)=0$ and such that
$$
F(t,y,r,\a\exp_{y}^{-1}(x),Q)-F(t,x,r,-\a\exp_{x}^{-1}(y),P)\leq
\omega(\a d^{2}(x,y)+d(x,y)),\eqno{(3.1)}
$$
for all fixed $t\in (0,T)$ and for all $x, y\in \o, r\in R, P\in
{\cal{L}}^{2}_{s}(TM_{x}), Q\in {\cal{L}}^{2}_{s}(TM_{y})$ with
$$
-(\frac{1}{\ep_{\a}}+\|A_{\a}\|)\left(
\begin{array}{cc}
I & 0\\
0 & I
\end{array}
\right)\leq \left(
\begin{array}{cc}
P & 0\\
0 & -Q
\end{array}
\right)\leq A_{\a}+\ep_{\a}A_{\a}^{2}, \eqno{(3.2)}
$$
where $A_{\a}$ is the second derivative of the function
$\p_{\a}(x,y)=\frac{\a}{2}d^{2}(x,y)$ at the point $(x,y)\in M\times
M$,
$$
\ep_{\a}=\frac{1}{2(1+\|A_{\a}\|)}
$$
and the points $x, y$ are assumed to be close enough to each other
so that  $d(x,y)<\min\{i_{M}(x),i_{M}(y)\}$.

Let $u\in USC([0,T)\times \bar{\o})$ be a subsolution and $v\in
LSC([0,T)\times \bar{\o})$ a supersolution of (1.1). Then $u\leq v$
on $[0,T)\times \o$.

In particular the Cauchy-Dirichlet problem (1.1) has at most one
viscosity solution.

{\bf Proof. }We first observe that for $\ep>0$,
$\tilde{u}=u-\frac{\ep}{T-t}$ is a vicosity subsolution of the
problem with the form
$$\left\{
\begin{array}{l}
\tilde{u}_{t}+F(t,x,\tilde{u}+\frac{\ep}{T-t},d\tilde{u},d^{2}\tilde{u})=-\frac{\ep}{(T-t)^{2}} \mbox{ \ in \ }(0,T)\times \o,\\
\tilde{u}(t,x)=h(t,x)+\frac{\ep}{t-T}, (t,x)\in [0,T)\times \partial \o, \\
 u(0,x)=\psi(x)-\frac{\ep}{T},  x\in \bar{\o},\\
\lim\limits_{t\uparrow T}\tilde{u}(t,x)=-\infty \mbox{ \ uniformly
on \ } \bar{\o}.
\end{array}
\right.
$$
Since $u\leq v$ follows from $\tilde{u}\leq v$ in the limit
$\ep\downarrow 0$, it will simply suffice to prove that
$$\tilde{u}(t,x)\leq v(t,x), (t,x)\in (0,T)\times \o.$$

We will assume $$\exists (s,z)\in (0,T)\times \o \mbox{ \ and \ }
\tilde{u}(s,z)-v(s,z)=\d>0$$ and then contradict this assumption.

By compactness of $[0,T]\times \bar{\o}\times \bar{\o}$ and upper
semicontinuity of $u-v$, and considering $-\frac{\ep}{T-t}\downarrow
-\infty$ when $t\uparrow T$, we have $\tilde {u}-v$ is bounded above
on $[0,T)\times \bar{\o}\times \bar{\o}$. Thus there exists
$(t_{\a},x_{\a},y_{\a})\in [0,T)\times \bar{\o}\times \bar{\o}$ such
that it is a maximum point of
$\tilde{u}(t,x)-v(t,y)-\frac{\a}{2}d^{2}(x,y)$ on $[0,T)\times
\bar{\o}\times \bar{\o}$. We set
$$m_{\a}=\tilde{u}(t_{\a},x_{\a})-v(t_{\a},y_{\a})-\frac{\a}{2}d^{2}(x_{\a},y_{\a}),$$
and obviously
$$
m_{\a}\geq \tilde{u}(s,z)-v(s,z)-\frac{\a}{2}d^{2}(z,z)=\d>0.
$$

Let us admit for a moment the following lemma:

\vskip 3mm

 {\bf Lemma 3.2. }

(i)$t_{\a}\neq 0,$

(ii)There exists $(t_{0},x_{0}) \in [0,T]\times \o$ s.t. when
$\a\rightarrow +\infty$, along some subsequence which we will still
denote by $(t_{\a},x_{\a},y_{\a})$ s.t.
$(t_{\a},x_{\a},y_{\a})\rightarrow (t_{0},x_{0},x_{0})$ and $\a
d^{2}(x_{\a},y_{\a})\rightarrow 0$.

\vskip 3mm

By Lemma 3.2, we have $t_{\a}\in (0,T),x_{\a},y_{\a}\in \o$ for
large $\a$.

According to [2], there is $r_{0}>0$ such that for every $x,y\in
B(x_{0},r_{0})$, we have that $d(x,y)<\min\{i_{M}(x),i_{M}(y)\}$,
the vectors $\exp_{x}^{-1}(y)\in TM_{x}\equiv TM_{x}^{*}$ and
$\exp_{y}^{-1}(x)\in TM_{y}\equiv TM_{y}^{*}$ are well defined, and
the function $d^{2}(x,y)$ is $C^{2}$ smooth on $B(x_{0},r_{0})\times
B(x_{0},r_{0})\in M\times M$. And we can assume that
$x_{\a},y_{\a}\in B(x_{0},r_{0})$ for all $\a$.

Now, for each $\a$, we can apply Theorem 2.7 with
$\o_{1}=\o_{2}=B(x_{0},r_{0}), u_{1}=\tilde{u}, u_{2}=-v,
\p(t,x,y)=\p_{\a}(x,y)=\frac{\a}{2}d^{2}(x,y)$, and for
$$
\ep=\ep_{\a}:=\frac{1}{2(1+\|d^{2}\p_{\a}(x_{\a},y_{\a})\|)}.
$$
Since $(t_{\a},x_{\a},y_{\a})$ is a global maximum of the function
$\tilde{u}(t,x)-v(t,y)-\frac{\a}{2}d^{2}(x,y)$ on $(0,T)\times
\o_{1}\times \o_{2}$, the condition $(2.1)$ in Theorem 2.7 is
guaranteed by having $\tilde{u}$ (and $v$) be a subsolution (resp.
supersolution) of a parabolic equation, thus for $\ep=\ep_{\a}$
there are numbers $b_{1},b_{2}$ and bilinear forms $P\in
{\cal{L}}^{2}_{s}(TM_{x_{\a}}), Q\in {\cal{L}}^{2}_{s}(TM_{y_{\a}})$
s.t.
$$
(b_{1},\frac{\partial}{\partial x}\p(t_{\a},x_{\a},y_{\a}),P)\in
\bar{{\cal{P}}}^{2,+}\tilde{u}(t_{\a},x_{\a}),
(-b_{2},-\frac{\partial}{\partial y}\p(t_{\a},x_{\a},y_{\a}),Q)\in
\bar{{\cal{P}}}^{2,-}v(t_{\a},y_{\a}),
$$
and
$$
b_{1}+b_{2}=\p_{t}(t_{\a},x_{\a},y_{\a})=0,
-(\frac{1}{\ep_{\a}}+\|A_{\a}\|)\left(
\begin{array}{cc}
I & 0\\
0 & I
\end{array}
\right)\leq \left(
\begin{array}{cc}
P & 0\\
0 & -Q
\end{array}
\right)\leq A_{\a}+\ep_{\a}A_{\a}^{2},
$$
where $A_{\a}=d^{2}\p(t_{\a},x_{\a},y_{\a})\in
{\cal{L}}^{2}_{s}(T(M\times M)_{(x_{\a},y_{\a})})$. Therefore,
according to the assumption, we have that
$$
F(t_{\a},y_{\a},r,\a\exp_{y_{\a}}^{-1}(x_{\a}),Q)-F(t_{\a},x_{\a},r,-\a\exp_{x_{\a}}^{-1}(y_{\a}),P)\leq
\omega(\a d^{2}(x_{\a},y_{\a})+d(x_{\a},y_{\a})). \eqno{(3.3)}
$$
On the other hand, by the conclusion in Section 3 of [2], we have
that
$$
\frac{\partial}{\partial
x}\p(t_{\a},x_{\a},y_{\a})=-\a\exp_{x_{\a}}^{-1}(y_{\a}),
-\frac{\partial}{\partial
y}\p(t_{\a},x_{\a},y_{\a})=\a\exp_{y_{\a}}^{-1}(x_{\a}).
$$
Combining the definition of viscosity supersolution and subsolution,
we have
$$
\begin{array}{l}
b_{1}+F(t_{\a},x_{\a},\tilde{u}(t_{\a},x_{\a})+\frac{\ep}{T-t},-\a\exp_{x_{\a}}^{-1}(y_{\a}),P)\leq
-\frac{\ep}{(T-t_{\a})^{2}}\leq -\frac{\ep}{T^{2}}<0,\\
-b_{2}+F(t_{\a},y_{\a},v(t_{\a},y_{\a}),\a\exp_{y_{\a}}^{-1}(x_{\a}),Q)\geq
0,
\end{array}
$$
this with properness of $F$ and (3.3), we have
$$
\begin{array}{ll}
\frac{\ep}{T^{2}}& \leq
F(t_{\a},y_{\a},v(t_{\a},y_{\a}),\a\exp_{y_{\a}}^{-1}(x_{\a}),Q)
-F(t_{\a},x_{\a},\tilde{u}(t_{\a},x_{\a})+\frac{\ep}{T-t},-\a\exp_{x_{\a}}^{-1}(y_{\a}),P)\\
& \leq
F(t_{\a},y_{\a},\tilde{u}(t_{\a},x_{\a}),\a\exp_{y_{\a}}^{-1}(x_{\a}),Q)
-F(t_{\a},x_{\a},\tilde{u}(t_{\a},x_{\a}),-\a\exp_{x_{\a}}^{-1}(y_{\a}),P)\\
& \leq \omega(\a d^{2}(x_{\a},y_{\a})+d(x_{\a},y_{\a})).
\end{array}
$$
This is a contradiction by Lemma 2.2 if we let $\a$ tend to
$\infty$.

\vskip 3mm
 The proof of Lemma 3.2:

(i)If $t_{\a}=0$, we have
$$
0<\d \leq
m_{\a}=\psi(x_{\a})-\psi(y_{\a})-\frac{\a}{2}d^{2}(x_{\a},y_{\a})-\frac{\ep}{T}
\leq
\sup\limits_{\bar{\o}\times\bar{\o}}(\psi(x)-\psi(y)-\frac{\a}{2}d^{2}(x,y))-\frac{\ep}{T}.
$$
However, since $\psi\in C(\bar{\o})$, the right-hand side above
tends to $-\frac{\ep}{T}$ according to Lemma 4.1 in [1], so
$t_{\a}\neq 0$ if $\a$ is large.

(ii)Since $(t_{\a},x_{\a},y_{\a})$ is a sequence in a compact set
$[0,T]\times \bar{\o}\times \bar{\o}$, there must be some
subsequence which we still denote $(t_{\a},x_{\a},y_{\a})$ that is
convergent to a limit $(t_{0},x_{0},y_{0})\in [0,T]\times
\bar{\o}\times \bar{\o}$ when $\a\rightarrow +\infty$. If $x_{0}\neq
y_{0}$, $\a d^{2}(x_{\a},y_{\a})$ will tends to $+\infty$. Combining
that $\tilde{u}-v$ is upper bounded, we get $m_{\a}\rightarrow
-\infty$, which is also a contradiction to $m_{\a}\geq \d>0$.

 When $t_{\a}\uparrow T , m_{\a}\rightarrow
-\infty$, so $t_{0}\neq T$. If $x_{0}=y_{0}\in \partial \o$, by the
upper semicontinuity of $\tilde{u}-v$, we have
$$
\begin{array}{ll}
0<\d \leq \overline{\lim\limits_{\a\rightarrow
\infty}}m_{\a}&=\overline{\lim\limits_{\a\rightarrow
\infty}}(\tilde{u}(t_{\a},x_{\a})-v(t_{\a},y_{\a})-\frac{\a}{2}d^{2}(x_{\a},y_{\a}))\\
& \leq \overline{\lim\limits_{\a\rightarrow
\infty}}(\tilde{u}(t_{\a},x_{\a})-v(t_{\a},y_{\a}))\\
& \leq \tilde{u}(t_{0},x_{0})-v(t_{0},x_{0})\leq
\frac{\ep}{t-T}-0<0,
\end{array}
$$
which is obviously a contradiction, so $x_{0}\in \o$.

Since $\tilde{u}-v$ is bounded above, there exists a constant $C$
s.t.
$$
0<\d\leq
m_{\a}=\tilde{u}(t_{\a},x_{\a})-v(t_{\a},y_{\a})-\frac{\a}{2}d^{2}(x_{\a},y_{\a})
\leq \tilde{u}(t_{\a},x_{\a})-v(t_{\a},y_{\a}) \leq C.
$$
So there exists a subsequence such that $m_{\a}$ converges.
According to Proposition 3.7 in [1], we have $\a
d^{2}(x_{\a},y_{\a})\rightarrow 0.$

{\bf Remark 3.3. }(see [2])If $M$ has nonnegative sectional
curvature, then condition $(3.2)$ implies that $P\leq L_{yx}(Q)$.

Therefore, if $M$ has nonnegative curvature and $F$ is degenerate
elliptic for each fixed $t$, then $(3.2)$ automatically implies that
$$
F(t,x,r,-\a\exp^{-1}_{x}(y),L_{yx}Q)-F(t,x,r,-\a\exp^{-1}_{x}(y),P)\leq
0.
$$
hence
$$
\begin{array}{ll}
& F(t,y,r,\a\exp^{-1}_{y}(x),Q)-F(t,x,r,-\a\exp^{-1}_{x}(y),P)\\
= & F(t,y,r,\a\exp^{-1}_{y}(x),Q)-F(t,x,r,-\a\exp^{-1}_{x}(y),L_{yx}Q)\\
& +F(t,x,r,-\a\exp^{-1}_{x}(y),L_{yx}Q)-F(t,x,r,-\a\exp^{-1}_{x}(y),P)\\
\leq &
F(t,y,r,\a\exp^{-1}_{y}(x),Q)-F(t,x,r,-\a\exp^{-1}_{x}(y),L_{yx}Q),
\end{array}
$$
and we see that the main condition $(3.1)$ in Theorem 3.1 is
satisfied if we automatically require, for instance, that
$$
F(t,y,r,\a\exp^{-1}_{y}(x),Q)-F(t,x,r,-\a\exp^{-1}_{x}(y),L_{yx}Q)\leq
\omega(\a d^{2}(x,y)+d(x,y)).
$$
So we only need that $F$ is $\emph{intrinsically uniformly
continuous}$ w.r.t. $x$ uniformly in $t$(see [2]).

{\bf Definition 3.4. }We will say that $F\in C(\chi,R)$ is
$\emph{intrinsically uniformly continuous}$ w.r.t. $x$ uniformly in
$t$, if
$$
\sup_{t\in
(0,T)}|F(t,y,r,L_{xy}\z,L_{xy}P)-F(t,x,r,\z,P)|\rightarrow 0 \mbox{
\ uniformly as \ }y\rightarrow x.
$$

 Let us sum up what we have just shown:

\vskip 5mm

{\bf Corollary 3.5. }Let $\o$ be a bounded open subset of a complete
finite-dimensional Riemannian manifold $M$ with nonnegative
sectional curvature, and $F \in C(\chi, R)$ be continuous, proper
for each fixed $t\in (0,T)$ and intrinsically uniformly continuous
w.r.t. $x$ uniformly in $t$.

Let $u\in USC([0,T)\times \bar{\o})$ be a subsolution and $v\in
LSC([0,T)\times \bar{\o})$ a supersolution of (1.1). Then $u\leq v$
on $[0,T)\times \o$.

In particular PDEs (1.1) has at most one viscosity solution.

\vskip 5mm

When $M$ has negative curvature, the main condition (3.1) in Theorem
3.1 involves kind of a uniform continuity assumption on the
dependence of $F$ with respect to $d^{2}u(t,x)$ by virtue of the
following remark(see [2]).

{\bf Remark 3.6. }Assume that $M$ has sectional curvature bounded
below by some constant $-K_{0}\leq 0$. Then condition (3.2) in
Theorem 3.1 implies that
$$
P-L_{yx}(Q)\leq \frac{3}{2}K_{0}\a d^{2}(x,y)I,
$$
where $I(v,v)=\|v\|^{2}$.

 {\bf Corollary 3.7. }Let $M$ be a compact
Riemannian manifold (no assumption on curvature) and $\o$ is a
bounded open subset of $M$. Suppose that $F\in C(\chi,R)$ is proper,
continuous, and satisfies the following uniform continuity
assumption: for every $\ep>0$, there exists $\d>0$ such that
$$
d(x,y)\leq \d, P-L_{yx}(Q)\leq \d I \Rightarrow \sup_{t\in
(0,T)}\{F(t,y,r,\a\exp^{-1}_{y}(x),Q)-F(t,x,r,-\a\exp^{-1}_{x}(y),P)\}\leq
\ep,
$$
for each fixed $t$ and for all $x,y\in \o$ with $d(x,y)<i(M), r\in
R, \a>0, P\in {\cal{L}}^{2}_{s}(TM_{x})$ and $Q\in
{\cal{L}}^{2}_{s}(TM_{y})$. Then there is at most one viscosity
solution to PDEs (1.1).

Remind that $M$ has sectional curvature bounded below since it is
compact, so the conclusion is an analogue of Corollary 4.10 in [2].
We omit the proof.

\vskip 1cm

\section{Existence result}\label{sec:intro}

\qquad In [1], detailed research called Perron's method has been
adapted to establish existence of viscosity solution to the
Dirichlet problem in Euclidean space. The same method has been used
to get the existence result for Dirichlet problem on Riemannian
manifolds. For our PDEs (1.1), we can go exactly as in [1] with
appropriate changes to get our result as follows:

{\bf Theorem 4.1.} Let comparison hold for (1.1) , i.e., if $\omega$
is a subsolution of (1.1) and $v$ is a supersolution of (1.1), then
$\omega \leq v$. Suppose also that there exists a subsolution
$\b{u}$ and a supersolution $\bar{u}$ of (1.1) that satisfy the
initial condition $\b{u}_{*}(0,x)=\bar{u}^{*}(0,x)=\psi (x)$ for
$x\in \bar{\o}$ and $\b{u}_{*}(t,x)=\bar{u}^{*}(t,x)=h(t,x)$ for
$(t,x)\in [0,T)\times \partial{\o}$ . Then
$$
W(t,x)=\sup\{\omega(t,x) : \b{u}\leq \omega\leq \bar{u} \mbox{ \ and
\ } \omega \mbox{ \ is a subsolution of (1.1)} \}
$$
is a solution of (1.1).

Here we used the following notation:
$$
\begin{array}{l}
u^{*}(t,x)=\lim_{r\downarrow 0}\sup\{u(s,y) : (s,y)\in
(0,T)\times M \mbox{ \ and \ } |s-t|\leq r, d(y,x)\leq r\};\\
u_{*}(t,x)=\lim_{r\downarrow 0}\inf\{u(s,y) : (s,y)\in (0,T)\times M
\mbox{ \ and \ } |s-t|\leq r, d(y,x)\leq r\},
\end{array}
$$
that is $u^{*}$ denotes the upper semicontinuous envelope of $u$
(the smallest upper semicontinuous function, with values in
$[-\infty,+\infty]$, satisfying $u\leq u^{*}$), and similarly
$u_{*}$ stands for the lower semicontinuous envelope of $u$.

As in [1], we need the following several steps which are Riemannian
versions for parabolic PDEs.

{\bf Proposition 4.2.} Let $\o\subset M$ be locally compact, $u\in
USC([0,T)\times \bar{\o})$, $(t,z)\in (0,T)\times \o$, and $(p,\z,
A)\in {\cal{P}}^{2,+}u(t,z)$. Suppose that $u_{n}$ is a sequence of
upper semicontinuous functions defined on $[0,T)\times \bar{\o}$
s.t.
$$\begin{array}{l} \mbox{(i)there
exists \ } (t_{n},x_{n})\in (0,T)\times \o\mbox{ \ such that \ }
(t_{n},x_{n},u_{n}(t_{n},x_{n}))\rightarrow (t,z,u(t,z)),\\
\mbox{(ii)if \ }(s_{n},y_{n})\in (0,T)\times \o$ and
$(s_{n},y_{n})\rightarrow (s,y)\in (0,T)\times \o,\\
\mbox{ \ \ \  then \ }
\limsup_{n\rightarrow\infty}u_{n}(s_{n},y_{n})\leq u(s,y).
\end{array}$$
Then there exist $(\hat{t}_{n},\hat{x}_{n})\in (0,T)\times \o$ and
$(p_{n},\z_{n},A_{n})\in {\cal{P}}^{2,+}u_{n}(t_{n},z_{n})$ such
that
$$(\hat{t}_{n},\hat{x}_{n},u_{n}(\hat{t}_{n},\hat{x}_{n}),p_{n},\z_{n},A_{n})\rightarrow
(t,z,u(t,z),p,\z,A).$$

{\bf Proof.} The Euclidean version of the above proposition for
elliptic functions was proved in [1]. Azagra et al. generalized it
to the setting of Riemannian manifolds. For Euclidean version of
parabolic functions, we can refer you to [4]. If we note Remark 2.3,
we can also generalized the result in [4] to get our proof for the
above Proposition with the same method as in [2]. So we omit the
proof.

Use the above proposition, we can soon get the following lemma just
like in [1]:

{\bf Lemma 4.3.} Let $\o\in M$ be locally compact and  $F\in
C(\chi,R)$ be continuous for each fixed $t\in (0,T)$.   ${\cal{F}}$
is a family of solution of $u_{t}+F\leq 0$ in $(0,T)\times \o$. Let
$\omega(t,x)=\sup\{u(t,x) : u\in {\cal{F}}\}$ and assume that
$\omega^{*}(t,x)<\infty$ for $(t,x)\in (0,T)\times \o$. Then
$\omega^{*}$ is a solution of $u_{t}+F\leq 0$ in $(0,T)\times \o$.

The following lemma is also important to the existence result.

{\bf Proposition 4.4.} Let $\o\subset M$ be open and $u$ be solution
of $u_{t}+F\leq 0$ in $(0,T)\times \o$. If $u_{*}$ fails to be a
supersolution at some point $(\hat{t},\hat{x})$, i.e., there exists
$(p,\z,A)\in {\cal{P}}^{2,-}u_{*}(\hat{t},\hat{x})$ for which
$p+F(\hat{t},\hat{x},u_{*}(\hat{t},\hat{x}),\z,A)<0$, then for any
small $\kappa>0$ there is a subsolution $u_{\kappa}$ of $u_{t}+F\leq
0$ in $(0,T)\times \o$ satisfying
$$
\left\{
\begin{array}{l}
u_{\kappa}(t,x)\geq u(t,x) \mbox{ \ and \ } \sup_{(0,T)\times
\o}(u_{\kappa}-u)>0,\\
u_{\kappa}(t,x)=u(t,x)\mbox{ \ for \ } (t,x)\in (0,T)\times
\o,|t-\hat{t}|+d(x,\hat{x})\geq \kappa
\end{array}
\right.
$$

{\bf Proof.} Set
$$\p_{\d,\g}(t,v)=u_{*}(\hat{t},\hat{x})+\d+p(t-\hat{t})+\langle
\z,v\rangle_{\hat{x}}+\frac{1}{2}\langle
Av,v\rangle_{\hat{x}}-\g\|v\|^{2}-\g(t-\hat{t}).
$$
Then, by continuity, $\exists r,\d,\g>0$ small enough s.t. in
$$B_{r}=\{(t,v) : |t-\hat{t}|+\|v\|<r\},$$ we have
$$
p-\g+F(t,\exp_{\hat{x}}(v),\p_{\d,\g}(t,v),\z-2\g v,A-2\g I)<0.
$$
Set $x=\exp_{\hat{x}}(v)$. Since $(p,\z,A)\in
{\cal{P}}^{2,-}u_{*}(\hat{t},\hat{x})$, by Proposition 2.2 we have
$$
u(t,x)\geq u_{*}(t,x)= u_{*}(t,\exp_{\hat{x}}(v))\geq
u_{*}(\hat{t},\hat{x})+p(t-\hat{t})+\langle
\z,v\rangle_{\hat{x}}+\frac{1}{2}\langle
Av,v\rangle_{\hat{x}}-o(\|v\|^{2}+|t-\hat{t}|),
$$
if we choose $\d=\frac{r^{2}}{8}\g$ then
$u(t,x)>\p_{\d,\g}(t,\exp_{\hat{x}}^{-1}(x)):=u_{\d,\g}(t,x)$ for
$\frac{r}{2}\leq \|v\|\leq r$ if $r$ is sufficiently small.

According to Proposition 2.8 in [2], we have
$$
d\p_{\d,\g}(t,0_{\hat{x}})=du_{\d,\g}(t,\hat{x}),d^{2}\p_{\d,\g}(t,0_{\hat{x}})=d^{2}u_{\d,\g}(t,\hat{x}),
$$
and
$$
\frac{\partial \p_{\d,\g}}{\partial t}=\frac{\partial
u_{\d,\g}}{\partial t},\lim_{v\rightarrow
0_{\hat{x}}}d\p_{\d,\g}(t,v)=\lim_{x\rightarrow
\hat{x}}du_{\d,\g}(t,x),\lim_{v\rightarrow
0_{\hat{x}}}d^{2}\p_{\d,\g}(t,v)=\lim_{x\rightarrow
\hat{x}}d^{2}u_{\d,\g}(t,x),
$$
in a neighborhood of $(\hat{t},\hat{x})$.

So there exists $\kappa>0$ small enough s.t. $u_{\d,\g}\in
C^{1,2}(B_{\kappa})$, where $B_{\kappa}=\{(t,x) :
|t-\hat{t}|+d(x,\hat{x})<\kappa\}$. Moreover,
$u(t,x)>u_{\d,\g}(t,x)$ and $u_{\d,\g}(t,x)$ is a solution of
$u_{t}+F\leq 0$ when $(t,x)\in B_{\kappa}$. Then, by Lemma 4.3, the
function
$$
u_{\kappa}(t,x)=\left\{\begin{array}{ll}
\max(u(t,x),u_{\d,\g}(t,x)),
&(t,x)\in B_{\kappa},\\
u(t,x)& \mbox{otherwise},
\end{array}
\right.
$$
is a solution of $u_{t}+F\leq 0$ in $(0,T)\times \o$. The last
observation is that in every neighborhood of $(\hat{t},\hat{x})$
there are points such that $u_{\kappa}(t,x)>u(t,x)$; indeed, by
definition, there is sequence $(t_{n},x_{n},u(t_{n},x_{n}))$
convergent to $(\hat{t},\hat{x},u_{*}(\hat{t},\hat{x}))$ and then
$$
\lim_{n\rightarrow
\infty}(u_{\kappa}(t_{n},x_{n})-u(t_{n},x_{n}))=u_{*}(\hat{t},\hat{x})+\d-u_{*}(\hat{t},\hat{x})>0.
$$
$$
\eqno{\Box}
$$

Having these above preparation, we can easily achieve our proof of
Theorem 4.1 as in [1]. We omit it.

\end{document}